\documentclass[12pt,a4paper]{article}
\usepackage[utf8]{inputenc}
\usepackage{amsmath}
\usepackage{amsfonts}
\usepackage{amssymb}
\newtheorem{theo}{Théorème}[section]
\newtheorem{cj}{Conjecture}[section]
\newtheorem{lm}{Lemme}[section]
\newtheorem{cor}{Corollaire}[section]
\usepackage{graphics,color}
\usepackage{graphicx}
\usepackage{xr-hyper} 
\usepackage{hyperref} 
\hypersetup{                    % parametrage des hyperliens
    colorlinks=true,                % colorise les liens
    breaklinks=true,                % permet les retours à la ligne pour les liens trop longs
    urlcolor= blue,                 % couleur des hyperliens
    linkcolor= blue,                % couleur des liens internes aux documents (index, figures, tableaux, equations,...)
    citecolor= green                % couleur des liens vers les references bibliographiques
    }

\begin{document}

\begin{center}
\textbf{MESURE ET ACTION DES I-PERMUTATIONS SUR LES \textcolor{red}{M}\textcolor{cyan}{U}\textcolor{yellow}{L}\textcolor{magenta}{T}\textcolor{black}{I}\textcolor{green}{G}\textcolor{red}{R}\textcolor{blue}{A}\textcolor{cyan}{P}\textcolor{magenta}{H}\textcolor{black}{E}\textcolor{blue}{S} MULTICOLORES FINIS ET INIFINIS}

\textcolor{red}{Version \footnote{\textcolor{blue}{ Les raisons d'introduction des dilatations : Section \ref{limite} } } du 1 Juillet 2015  }

Présentée à :

Université Claude-Bernard-Lyon1. Département de Mathématiques. 43, boulevard du 11 novembre 1918, F-69621-Villeurbanne, France.

Par : Mohamed Sghiar

 msghiar21@gmail.com
 
 Tel : 0033(0)953163155 \& 0033(0)669753590.

\end{center}

\textbf{Abstract  } : Among other results, the purpose of this article is to show the existence of an  $\mathbb{R}$-space-vector with basis $\omega^i_j $,  i, j are integers such that every graph  with n vertex   $ n \geq  3 $ is the vector: $$\mathcal{V}(n) = \sum_{j = 0}^{n-1} {\alpha^{n-1} _j} \omega ^ {n-1} _j $$ Where $ {\alpha ^ {n-1} _j }$ is the number of sub graphs of type $\omega^{n-1}_j$ . We deduce that two graphs are isomorphic if for any measure, they have the same number of maximal proper subset with this measure.

\textbf{Résumé} : Entre autres résultats, le but de cet article est montrer l'existence d'un $\mathbb{R}$-espace-vectoriel de base  $\omega^i_j $ où i, j sont des entiers, tel que tout graphe $\mathcal{V}$  de cardinal $ n \geq 3 $ est le vecteur : $$\mathcal{V}(n) = \sum_{j=0}^{n-1}{\alpha^{n-1}_j}\omega^{n-1}_j $$ Où  ${\alpha^{n-1}_j}$ est le nombre  de sous  graphes de type $\omega^{n-1}_j $. On en déduit que deux graphes sont isomorphes si pour toute mesure, ils ont le même nombre de parties propres maximales ayant cette mesure.

\newpage
\renewcommand{\contentsname}{Sommaire} % Dans le corps du document,avant la commande \tableofcontents.
\tableofcontents

\newpage

\begin{center}
\section*{Remerciements :}
\end{center}

\addcontentsline{toc}{section}{Remerciements} 

Je tiens à remercier le professeur Maurice Pouzet pour sa lecture, ses conseils et ses suggestions.

Je remercie aussi toute personne qui contribue à la réussite des résultats de cette   œuvre dont les techniques ont permis de résoudre  de nombreuses conjectures aussi bien pour le cas des graphes finis que  pour le cas des graphes  infinis.... des techniques qui s'appliquent donc de "l'infiniment petit à l'infiniment grand "...

\begin{flushright}
sghiar
\end{flushright}

\newpage

\begin{center}

\section*{Introduction}

\end{center}

\addcontentsline{toc}{section}{Introduction}

Dans le théorème \ref{t1} je démontre que  toute i-permutation sur E ( c.a.d une permutation sur les parties à i éléments de E) est déduite de l'action d'une permutation sur les éléments de E.

Puis, je donne dans le théorème \ref{t2} une application des i-permutations dans la preuve de la conjecture d'Ulam [1 et 10].

Après avoir introduit la notion de mesure sur les graphes,  je donne dans le théorème \ref{th mes}, une généralisation de la conjecture d'Ulam, à savoir que  deux graphes sont isomorphes dès que pour toute mesure, ils ont le même nombre de parties propres maximales ayant cette mesure.

Dans le corollaire \ref{c2} , je donne une représentation vectorielle des graphes à au moins 3 éléments : plus précisément  je démontre l'existence d'un $\mathbb{R}$-espace-vectoriel de base  $\omega^i_j $ où i, j sont des entiers, tel que tout graphe $\mathcal{V}$  de cardinal $ n \geq 3 $ est le vecteur $$\mathcal{V}(n) = \sum_{j=0}^{n-1}{\alpha^{n-1}_j}\omega^{n-1}_j $$ Où  ${\alpha^{n-1}_j}$ est le nombre  de sous  graphes de type $\omega^{n-1}_j $.

Et dans le corollaire \ref{c3} on déduit que toute (n-1)-permutation $ \sigma_{n-1}$ sur E  et préservant les mesures entre deux graphes sur E  est déduite (à une permutation près) d'une permutation $ \sigma $  sur E et préservant toutes les mesures entre ces deux graphes.

Dans la deuxième section je démontre que les principaux résultats de la premier section se généralisent aux graphes multicolores.

Dans la troisième section je réponds à une question de Maurice Pouzet sur la reconstruction des multigraphes.

Dans la quatrième section, je donne une réponse positive à la conjecture de M. Pouzet sur le reconstruction des relations m-aires h-symétriques.

Dans la cinquième section je démontre les raisons de la limite de validité de l'utilisation des $\infty$ -permutations pour les graphes infinies. 

Toutefois, j'énonce dans la sixième section certaines conjectures sur les multigraphes multicolores infinis. 

Dans les  septième et huitième sections, j'utilise la\textbf{ dilatation} des  $\infty$-permutations pour résoudre les conjectures ci-dessus, en particulier la reconstruction des graphes infinis et  localement finis.

\newpage

\begin{center}

\section*{Notations et  définitions}

\end{center}

\addcontentsline{toc}{section}{Notations et  définitions}

Soit $\omega^{i}_j,  $   le type d'un  graphe  de cardinal i (à un isomorphe près).

A tout $\omega^{i}_j,  $  associons un nombre $\mu ^{i}_{j}$, de tels  façon que $\mu ^{i}_{j} = \mu ^{k}_{l} \Longleftrightarrow i=k\,et \, j=l$

Pour tout graphe G de base E,  la mesure $\underset{G}{\mu} $  ou $\mu$  est la fonction  définie sur les parties de E par :  $\mu(\chi)=\mu ^{i}_{j}  $ si $G|\chi $  est isomorphe à $\omega^{i}_j$, 

$\mu ^{i}_{j}  $ est dite la mesure de $\chi$

Une arête {x, y} est dite pleine si G(x,y)=1, et elle est dite vide si G(x,y)= 0.

$\underset {G}{deg}(x)= card\{ y \in E \diagdown \{x\} \diagup G(x,y) = 1\}$

Une i-permutation $\sigma_i$ sur un ensemble E est une permutation sur les parties à i éléments de E.

Une i-permutation $\sigma_i$ sur les parties de E à i éléments est dite déduite d'une permutation $ \sigma $  sur E si: $ \sigma_i\chi = \sigma \chi , \; \forall\chi \subseteq E $ .

Si G et G' sont deux graphes sur un même ensemble E, une (n-1)-permutation entre les parties à $n-1$ éléments est dite préservant les mesures $\mu ^{n-1}_{i}$ si  $\underset{G}{\mu} (\chi)=\mu ^{n-1}_{i} \Longleftrightarrow \underset{G'}{\mu} (\sigma_ {n-1}(\chi))=\mu ^{n-1}_{i}$

Un graphe est dit bicolore si $ \text{card} \{\mu(\chi); \; \chi \subseteq E \diagup| \chi | = 2  \} \leq 2 $ et $ \text{card} \{\mu(\chi); \;\chi \subseteq E \diagup | \chi | = 1  \} =1 $.

Un graphe est dit multicolore si $ \text{card} \{\mu(\chi); \;\chi \subseteq E \diagup | \chi | = 2  \} \geq 2 $ et $ \text{card} \{\mu(\chi); \;\chi \subseteq E \diagup | \chi | = 1  \} =1 $.

Un multigraphe multicolore $(G_1 \cdots G_k)$ est un ensemble de k graphes multicolores sur une même base E.

Exemple de  Graphes et de multigraphes  multicolores : 
\begin{figure}
\centering
\includegraphics[width=50mm,height=25mm,scale=0.1]{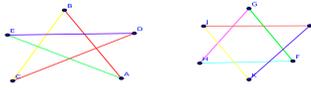}
\caption{Deux Graphes multicolores}

\end{figure}

\begin{figure}
\centering
\includegraphics[width=100mm,height=100mm,scale=0.1]{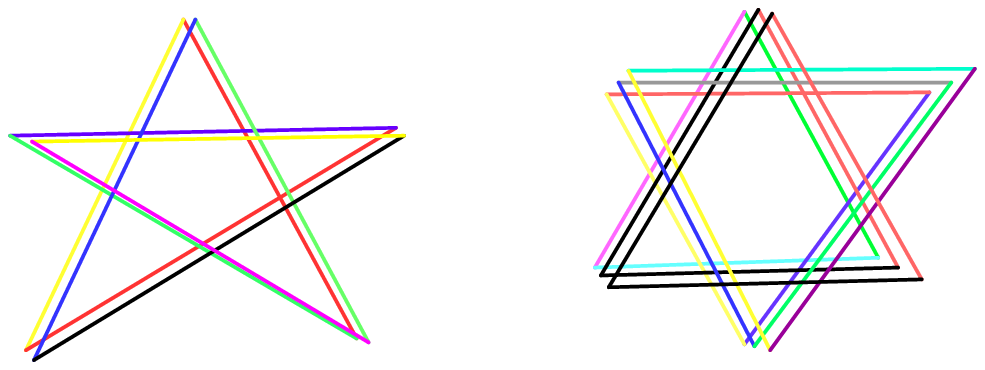}
\caption{Deux multigraphes multicolores}

\end{figure}

\newpage

\label{IF1}

\begin{center}
 \section*{ Une idée physique et géométrique}
 \end{center} 

\addcontentsline{toc}{section}{ Une idée physique et géométrique}

\textbf{Cas fini : }

Soient n particules $x_i, i \in \{1, \cdots, n\}$ , $(n \geq 3)$ .
Supposons que suite à une explosion les particules se changent de places en préservant les mesures $\mu^2_j$.

Dans ce cas, il est tout à fait normal que toutes les mesures $\mu^i_j$ resteront inchangées après l'explosion.

Si maintenant et inversement après  l'explosion, on a autant  de mesures  $\mu^{n-1}_j$ -sur les parties propres et maximales- qu'avant et qu'après l'explosion, on verra que pour le cas fini : l'explosion agit sur les particules par permutation tout en préservant toutes les mesures  $\mu^i_j$.

Autrement dit : Si G et G' sont les graphes représentant les particules avant et après l'explosion, alors G et G' sont isomorphes.

\textbf{Cas infini :}

Le cas infini est un peu plus délicat  mais il montre à quel point les techniques utilisées sont plausibles :

Si suite à une explosion les particules se changent de places en préservant les mesures $\mu^2_j$,  il est tout à fait normal que toutes les mesures $\mu^i_j$ resteront inchangées après l'explosion.

Si maintenant et inversement suite à l' explosion, chaque  partie propre et maximale de particules a tendance à se changer en une partie contenant une partie propre et maximale tout en préservant les  mesures -sur les parties contenant au moins une partie propre et maximale- , l'explosion n' agit-elle pas sur les particules par permutation tout en préservant toutes les mesures  $\mu^i_j$ ?

Autrement dit :  Les forces partielles qui agissent sur les parties propres et maximales ne sont-elles pas dues à une unique force qui agit sur les particules par permutation tout en préservant les mesures ?

Je réponds à cette question en montrant dans quels cas la réponse est positive en écartant certains types de garphes comme ceux trouvés dans les contres-exemples de J. Fisher.

Les contres-exemples de J. Fisher  confirment  bien  que mes nouvelles techniques sont efficaces et plausibles  pour  la reconstruction des graphes  finis et infinis.  Et ce,  rien  qu'en pensant  que les deux graphes G et G' ne sont que les états des particules  avant et après une explosion qui laisse invariantes les mesures  sur les parties propres et maximales !

L'idée principale consiste donc à tenir compte de l'invariance des mesures et  du fait que sur chaque point -ou particule- agit  une force résultante  des autres forces partielles  agissantes sur les parties propres et maximales -considérées comme des parois- et  contenant le dit-point (ou particule).

Et vu que les forces partielles sont définies à une permutation près sur les parties propres et maximales (en préservant les mesures), la force résultante, elle aussi, n'est  définie qu'à une permutation  près sur les particules mais en préservant les mesures.

\newpage

\label{IF2}

\begin{center}

\section*{Vision dynamique de le reconstruction des graphes} 
\end{center}

\addcontentsline{toc}{section}{Vision dynamique de le reconstruction des graphes} 

Si G est un graphe sur une base E de cardinal $n \geq 3$, on considère les points de E comme des particules, et la restriction de G à toute partie F d'éléments  E de cardinal i comme une mesure  $\mu^i_j$  sur F -j dépend de F-.

Si G' est un autre graphe sur E, et si $\sigma $ est un isomorphisme entre G et G' sur E, alors il existe une permutation $\sigma^*$ (ou $\sigma_{n-1}$  si E est de cardinal n)  sur les parties propres et maximales et préservant les mesures.

Plus tard, on va définir la \textbf{dilatation} $\hat{\sigma^{*}}$  de $\sigma^{*}$ définie de  $\Omega^*_E$ sur  $\Omega_E$ par : $\hat{\sigma^{*}}(\chi)= \sigma^{*}(\chi )$ si $\sigma^{*}(\chi) \neq E$  dans  $\Omega_E/ G'$ .  Et $\hat{\sigma^{*}}(\chi)= E$ si $G'|E \simeq G'|\sigma^{*} (\chi)$ dans $\Omega_E/ G'$.  où on a noté : 

  $\Omega^*_E$ l'ensemble des parties propres et maximales de E . 
 
  $\Omega_E$ l'ensemble des parties de E abritant au moins une partie propre et maximale de E.

Si E est de cardinal fini, alors   $\sigma^*=\hat{\sigma^{*}}$

\textbf{Le problème inverse }: Si il existe -à une permutation près-  une action $\hat{\sigma^{*}}$  qui transforme une  partie propre et maximale en une autre partie contenant une partie propre et maximale et en préservant  les mesures $\mu^2_i$. Peut-on trouver une permutation $\sigma $  sur E qui conserve les mesures $\mu^2_i$ entre G et G'?

On verra dans quels cas la permutation $\sigma $ existe, en particulier si : $\sigma^*=\hat{\sigma^{*}}$.

\textbf{Vision dynamique de le reconstruction des graphes :}

Dans les cas où on a   : $\sigma^*=\hat{\sigma^{*}}$.

Les  particules $x_i$, sous les actions  $\sigma^*(E\diagdown \{x_j\})$,  vont se positionner aux points $\sigma(x_i)$ -à une permutation près- de sorte que les mesures $\mu^2_i$ entre G et G' restent conservées.

On peut donc voir $\sigma$  comme l'action  résultante  des actions $\sigma^*(E\diagdown \{x_j\})$  où $x_j \in E$.

\newpage

\section{Cas où les graphes  sont bicolores}

\subsection{Action des i-permutations sur un ensemble E }

\begin{theo} \label{t1}

Soit E un ensemble de cardinal $ n \geq 3 $ .

Toute $(n-1)$-permutation $\sigma_{n-1}$ sur les parties de E à $n-1$ éléments est déduite d'une permutation $ \sigma $  des éléments de E : c'est à dire : 

$ \sigma_{n-1}\chi = \sigma \chi , \; \forall\chi \subseteq E \;\;  \diagup   \mid \chi \mid =  n-1 $.

Et  $ \sigma $ vérifie : $\sigma(A \cap B) =  \sigma_{n-1}(A) \cap \sigma_{n-1}(B) $  si $\mid A \mid=\mid B \mid= n-1$ .

On en déduit que : $$ \sigma(x_i) =  \bigcap_{j \in \{1,\cdots ,n\} \diagdown \{i\}}^{} {\sigma_{n-1}(E\diagdown \{x_j\})} \;  \forall x_i \in E $$

\end{theo}

\textbf{Preuve : }

\textbf{Première démonstration :
}
Supposons que le résultat est vrai jusqu'au  n-1 , et montrons qu'il est vrai pour n.
 
Soit x un élément de E.
Notons $C^i_E$  l'ensemble des parties de E à i éléments  et $C^i_{E,x}$  l'ensemble des parties de E à i éléments  contenant x.

On a :  $C^i_E = (C^i_{E,x})^{c}\cup C^i_{E,x}$ .

Posons $ Y= \sigma_i((C^i_{E,x})^{c})$.

Si $\forall y  \in E $  il existe au moins une partie $X_y$  de $(C^i_{E,x})^{c}$  telle que $y \in \sigma_i(X_y) $, alors : Comme $Y=\underset {z \in E} { \cup}(C^i_{E,z}\cap Y)$ , on a alors : $Y=\underset {z \in E} { \cup}( C^i_{E,z} \diagdown \sigma_i(C^i_{E,x})) $ et $Y=(\underset {z \in E} { \cup} C^i_{E,z}) \diagdown \sigma_i(C^i_{E,x}) $, et par suite  : $ \mid Y \mid = \mid \underset {z \in E} { \cup} C^i_{E,z}\mid - \mid \sigma_i(C^i_{E,x})  \mid $

Or $\mid Y \mid = \mid \sigma_i((C^i_{E,x})^{c})\mid = \mid (C^i_{E,x})^{c}\mid = C^i_{n-1}$,  $ \mid \sigma_i(C^i_{E,x})  \mid = \mid C^i_{E,x}  \mid  = C^{i-1}_{n-1}$, et $\mid \underset {z \in E} { \cup} C^i_{E,z}\mid = nC^{i-1}_{n-1} $ .

Donc finalement on doit avoir : $ C^i_{n-1} = nC^{i-1}_{n-1} - C^{i-1}_{n-1} = (n-1)C^{i-1}_{n-1}$

Ce qui est impossible car $C^i_{n-1} =\frac{1}{i(n-1-i)} C^{i-1}_{n-1}$ si $ i \neq n-1$, et si $ i = n-1$ on doit avoir $1=(n-1)C^{i-1}_{n-1}$, ce qui est impossible.

On en déduit qu'il existe un élément y de E  tel que pour toute partie $X_y$  de $(C^i_{E,x})^{c}$,  $y \notin \sigma_i(X_y) $, soit $y \in \sigma_i(X_z)$  avec $X_z \in C^i_{E,x}$, et par suite $C^i_{E,y}= \sigma_i C^i_{E,x}$

Et $\sigma_i$  est aussi une transformation bijective  de $C^i_{E\diagdown \{x\}}$  sur $C^i_{E\diagdown \{y\}}$, qui induit une transformation bijective $ \sigma_{i-1} $  de $ C^{i-1}_{E\diagdown \{x\}}$  sur $C^{i-1}_{E\diagdown \{y\}}$

Si $i= n$, alors par récurrence, il existe une transformation bijective $ \sigma_{x,y} $ de  $E\diagdown \{x\}$  sur $E\diagdown \{y\}$ telle que $ \sigma_{i-1}$ est déduite de $ \sigma_{x,y} $ .

En prolongeant $ \sigma_{x,y} $  sur E par  $ \sigma(x)=y $ , on déduit que $\sigma_i$  est déduite de la transformation  $\sigma $ sur les éléments de E. $ \checkmark$

\textbf{Deuxième démonstration :
}

Si $\sigma_{n-1}$  est une (n-1)-permutation sur E, associons à  $\sigma_{n-1}$  la permutation $\widetilde{\sigma}_{n-1} $ sur E définie par : $\widetilde{\sigma}_{n-1}(x)= E \diagdown \sigma_{n-1}(E \diagdown \{x\}) \; \forall x \in E $.

On en déduit que $E \diagdown \widetilde{\sigma}_{n-1}(x)=  \sigma_{n-1}(E \diagdown \{x\})\; \forall x \in E $, et par suite  pour toute permutation $\sigma$ :

$\sigma^{-1} \sigma_{n-1}(E \diagdown \{x\})=\sigma^{-1}(E \diagdown \widetilde{\sigma}_{n-1}(x))=  E \diagdown \sigma^{-1}\widetilde{\sigma}_{n-1}(x) \; \forall x \in E $.

Or il existe $\sigma$ tel que $\sigma^{-1}\widetilde{\sigma}_{n-1}= 1_E$ ($1_E $ est l'élément neutre du groupe des permutations).

$\Longrightarrow \sigma^{-1} \sigma_{n-1}(E \diagdown \{x\})= E \diagdown \{x\} \; \forall x \in E$

$\Longrightarrow  \sigma_{n-1}(E \diagdown \{x\})= \sigma(E \diagdown \{x\}) \; \forall x \in E$

$\Longrightarrow \sigma_{n-1}\chi = \sigma \chi , \; \forall\chi \subseteq E \;\;  \diagup   \mid \chi \mid =  n-1 $.

D'où le résultat. $ \checkmark$

\textbf{Troisième démonstration :
}

Posons : $$ \widetilde {\sigma}_{n-1}(x_j)= E \diagdown \sigma_{n-1}(E\diagdown \{x_j\})$$

on a  $$   \bigcap_{j \in \{1,\cdots ,n\} \diagdown \{i\}}^{} {\sigma_{n-1}(E\diagdown \{x_j\})} = \bigcap_{j \in \{1,\cdots ,n\} \diagdown \{i\}}^{} (E\diagdown \widetilde {\sigma}_{n-1}(x_j)) = E \diagdown  \bigcup_{j \in \{1,\cdots ,n\} \diagdown \{i\}}^{}\widetilde{\sigma}_{n-1}(x_j)= \widetilde{\sigma}_{n-1}(x_i) $$

Montrons que $\widetilde {\sigma}$ prolonge $\sigma_{n-1}$  ou autrement dit que  $\sigma_{n-1}$ est déduite de  $\widetilde {\sigma}$.

Si $ x_i \in E \diagdown \{x_j\}$, alors : $$\widetilde{\sigma}_{n-1}(x_i)= \bigcap_{j \in \{1,\cdots ,n\} \diagdown \{i\}}^{} {\sigma_{n-1}(E\diagdown \{x_j\})}\subseteq \sigma_{n-1}(E\diagdown \{x_j\}) $$

Et le résultat s'en déduit. $ \checkmark$

\textbf{Quatrième démonstration :
}

Sachant qu'on a n parties à $n-1$ éléments de E, alors le nombre des (n-1)-permutations sur E est $n!$.

Or chaque permutation $\sigma$ sur E induit une (n-1)-permutation  $\sigma_{n-1}$ sur E définie par : 

$$\sigma_{n-1}\chi=\sigma \chi ; \; \forall \chi \subseteq E \diagup  \; |\chi|=n-1 $$

Et comme on a $n!$ permutation sur E, donc ces permutations induisent les $n!$  (n-1)-permutations sur E.

Ainsi toute (n-1)-permutation sur E est déduite d'une unique permutation sur E.$\checkmark$

\newpage

 \subsection{ Utilisation des i-permutations dans la preuve de la conjecture d'Ulam }

\begin{theo}[Conjecture d'Ulam] \label{t2}

Soit G et G' deux graphes sur une même base E de cardinal au moins égal à 3.

Si G et G' sont (-1)-hypomorphes, alors G et G' sont isomorphes .

\end{theo}

\textbf{Preuve :
}
On aura besoin du lemme suivant:

\begin{lm}  \label{l2}

Soit G et G' deux graphes multicolores sur un même ensemble  E à au moins 3 éléments. 

Si $ \sigma_ {n-1} $  est une (n-1)-permutation entre les parties à $n-1$ éléments et préservant les mesures $\mu ^{2}_{i}$ sur les parties à n-1 éléments, alors $ \sigma_ {n-1} $  est déduite d'une permutation $\widetilde{\sigma}$  avec $\widetilde{\sigma}\tau$  préservant elle aussi les mesures $\mu ^{2}_{i}$ où $\tau$  est une permutation sur les éléments de E.

Et à une \textbf{représentation près  }des graphes G et G', $ \sigma_ {n-1} \tau $  sera déduite de  $\widetilde{\sigma}\tau$ avec $\widetilde{\sigma}\tau$  préservant les mesures.

\end{lm}

\textbf{Preuve : }

Posons : $\sigma_{n-1} \left(  \begin{array}{c}
x_{1,1}^{i} \\ 
\vdots \\ 
x_{n-1,1}^{i}
\end{array}\right) = \left( \begin{array}{c}
\sigma_{n-1}(x_{1,1}^{i})  \\ 
\vdots \\ 
\sigma_{n-1}(x_{n-1,1}^{i}) \\ 

\end{array} \right)   $ 

L'action de $\sigma_{n-1}$  sur la  partie à n-1 éléments $ \left(  \begin{array}{c}  
x_{1,1}^{i} \\ 
\vdots \\ 
x_{n-1,1}^{i}
\end{array}\right) $  préservant la mesure $\mu ^{2}_{i}$ c'est à dire : $ G(x_{k,1}^{i}, x_{l,1}^{i})= G'(\sigma_{n-1}^{}(x_{k,1}^{i}),\sigma_{n-1}^{}(x_{l,1}^{i})) ; \; \forall\{k,l,i\}$

En dessinant G et G' sur E (à une permutation  près sur les éléments de E), alors de  la troisième preuve du théorème \ref{t1} , on a  :

$$\widetilde{\sigma}(x_i)= \bigcap_{k \in \{1,\cdots ,n\} \diagdown \{i\}}^{} {\sigma_{n-1}(E\diagdown \{x_k\})} \forall x_i \in E $$

Donc :  $\forall (x_i ; x_j) \in E^2$ , 

\begin{eqnarray*}
\{\widetilde{\sigma}(x_i)\}  \cup \{\widetilde{\sigma}(x_j)\} &=& \bigcap_{k \in \{1,\cdots ,n\} \diagdown \{i\}}^{} {\sigma_{n-1}(E\diagdown \{x_k\})} \cup \bigcap_{k \in \{1,\cdots ,n\} \diagdown \{j\}}^{} {\sigma_{n-1}(E\diagdown \{x_k\})} \\
&=&  \bigcap_{k \in \{1,\cdots ,n\} \diagdown \{i,j\}  }^{} \sigma_{n-1}(E\diagdown \{x_k\})
\end{eqnarray*}

On a  $G'(\widetilde{\sigma}(x_i), \widetilde{\sigma}(x_j)) = G'(\sigma_{n-1}(x_i),\sigma_{n-1}(x_j)) $. Et comme  $\sigma_{n-1} $  préserve les mesures $\mu ^{2}_{i}$ , alors :

$$ G(x_i, x_j)=  G'(\sigma_{n-1}(x_i),\sigma_{n-1}(x_j))$$

\begin{center}
(Autrement dit $\sigma_{n-1}$ pousse les $x_i$  en préservant les mesures )
\end{center}

Donc  : $$ G(x_i, x_j)=  G'(\widetilde{\sigma}(x_i),\widetilde{\sigma}(x_j))$$ (à une permutation  près sur les éléments de E).

Soit donc : $$ G(x_i, x_j)= G'(\widetilde{\sigma}\tau(x_i),\widetilde{\sigma}\tau(x_j)) , \; \forall (x_i , x_j) \in E^2$$  Où $\tau$ est une permutation.

Montrons maintenant que  à une \textbf{représentation près  }des graphes G et G', $ \sigma_ {n-1}\tau $  sera déduite de  $\widetilde{\sigma}\tau$.

En effet : 

$$ \sigma_{n-1}\chi = \widetilde{\sigma} \chi , \; \forall\chi \subseteq E \;\;  \diagup   \mid \chi \mid =  n-1 $$

On en déduit  que :

$$ \sigma_{n-1}\tau \tau^{-1}\chi = \widetilde{\sigma} \tau \tau^{-1}\chi , \; \forall\chi \subseteq E \;\;  \diagup   \mid \chi \mid =  n-1 $$

Et comme :

 $$ G(x_i, x_j)= G'(\widetilde{\sigma}\tau(x_i),\widetilde{\sigma}\tau(x_j)) ; \; \forall (x_i , x_j) \in E^2 $$

Alors :

$$ G(\tau^{-1} (x_i),  \tau^{-1}(x_j))= G'( \widetilde{\sigma}\tau\tau^{-1}(x_i),\widetilde{\sigma} \tau \tau^{-1}(x_j)); \; \forall (x_i , x_j) \in E^2$$

C'est à dire que $\sigma_{n-1}\tau$ est déduite de $\widetilde{\sigma} \tau$ avec $\widetilde{\sigma} \tau$  préservant les mesures $\mu ^{2}_{i}$  entre les graphes $\tau^{-1}G$  et $\tau^{-1}G'$. $\checkmark$

\textbf{Preuve du Théorème }

\textbf{Preuve } (Valable même pour les graphes multicolores) \label{2D1}

En utilisant le lemme \ref{l2} on  trouve  une permutation $ \sigma $ préservant les mesures $\mu ^{2}_{i}$, donc $G(x,z)= G'(\sigma(x), \sigma(z))$  pour tout couple d'éléments x et z  de  E. D'où le résultat. $ \checkmark$

--------------------------------

\newpage

\subsection{  Généralisation de la conjecture  d'Ulam }

Soit $\omega^{i}_j,  $   le type d'un  graphe  de cardinal i (à un isomorphe près).

A tout $\omega^{i}_j,  $  associons un nombre $\mu ^{i}_{j}$, de tels  façon que $\mu ^{i}_{j} = \mu ^{k}_{l} \Longleftrightarrow i=k\,et \, j=l$

Pour tout graphe G de base E,  la mesure $\underset{G}{\mu} $  ou $\mu$  est la fonction  définie sur les parties de E par :  $\mu(\chi)=\mu ^{i}_{j}  $ si $G|\chi $  est isomorphe à $\omega^{i}_j$, 

$\mu ^{i}_{j}  $ est dite la mesure de $\chi$

\begin{theo}\label{th mes}
Soit G et G' deux graphes sur une même base E de cardinal au moins égal à 3.

Si $ \forall j$  G et G' abritent le même nombre de parties de mesure $\mu^{n-1}_j $  alors G et G' sont isomorphes .

\end{theo}

\textbf{Preuve}

Des hypothèses il existe une (n-1)-permutation $\sigma_{n-1}$ préservant les mesures $\mu^{n-1}_i$. 
Il s'en suit que G et G' seront (-1)-hypomorphes, Et le résultat se déduit directement du Théorème \ref{t2}. $ \checkmark$

\newpage
\subsection{ Représentation vectorielle des graphes }

Soit  $\mathcal{E}$ un  $\mathbb{R}$-espace-vectoriel de base  $\omega^i_j $ où i, j sont des entiers.

En identifiant chaque $\omega^i_j $ à un graphe de type $\omega^i_j $ (à un isomorphe près), le théorème \ref{th mes}, permet d'écrire  tout  graphe $\mathcal{V}$  de cardinal $ n \geq 3 $ sous la forme du vecteur $$\mathcal{V}(n) = \sum_{j=0}^{n-1}{\alpha^{n-1}_j}\omega^{n-1}_j $$ Où   ${\alpha^{n-1}_j}$ est le nombre  de sous  graphes de type $\omega^{n-1}_j $.

\begin{cor} \label{c2}

Il existe un $\mathbb{R}$-espace-vectoriel de base  $\omega^i_j $ où i, j sont des entiers, tel que tout graphe $\mathcal{V}$  de cardinal $ n \geq 3 $ est le vecteur $$\mathcal{V}(n) = \sum_{j=0}^{n-1}{\alpha^{n-1}_j}\omega^{n-1}_j $$ Où  ${\alpha^{n-1}_j}$ est le nombre  de sous  graphes de type $\omega^{n-1}_j $.

\end{cor}

Ce corollaire \ref{c2} permet de retrouver le théorème \ref{th mes}

Le corollaire \ref{c3} suivant  va établir un lien entre la (n-1)-permutation $ \sigma_{n-1}$ sur E et la permutation établissant un isomorphisme entre G et G'.

\begin{cor} \label{c3}
Soit G et G' deux graphes sur une même base E de cardinal au moins égal à 3.

Toute (n-1)-permutation $ \sigma_{n-1}$ sur E  et préservant les mesures :  c'est à dire $\underset{G'}{\mu} (\sigma_ {n-1}(\chi))=\underset{G}{\mu} (\chi) $ pour toute partie $\chi $ à $n-1$ éléments de E est déduite  d'une permutation $ \sigma \tau^{-1} $  sur E avec $ \sigma $ préservant toutes les mesures : c'est à dire : $\underset{G'}{\mu} (\sigma(\chi))=\underset{G}{\mu} (\chi)  \forall\chi \subseteq E    $. Avec $ \sigma_{n-1} \tau \chi = \sigma  \chi , \; \forall\chi \subseteq E \;\;  \diagup   \mid \chi \mid =  n-1 $. où $\tau$ est une permutation sur E.

De plus  $ \sigma $ vérifie :

 $$ \sigma(x) =  \bigcap_{j \in \{1,\cdots, n\}}^{} {\sigma_{n-1}(E\diagdown \{\tau (x_j)\})} \; x_j\neq x \; \forall x \in E $$ 
 
Et du Théorème \ref{t1} en prolongeant $\sigma_ {n-1}$ par : $$\sigma_ {n-1}(x) =\bigcap_{j \in \{1,\cdots, n\}}^{} {\sigma_{n-1}(E\diagdown \{ x_j\})} \; x_j\neq x \; \forall x \in E  $$
On aura  $\underset{G'}{\mu} (\sigma_ {n-1}\tau \chi)=\underset{G}{\mu} (\chi) $ pour toute partie $\chi $ de E .

\end{cor}

\textbf{Preuve }: De ce qui précède, on a déjà montré  l'existence de $ \sigma $ préservant toute les mesures vu qu'on a montré dans le Théorème \ref{th mes} et dans le corollaire \ref{c2} que G et G' sont isomorphes. il reste  à démontrer que : $$ \sigma(x) =  \bigcap_{j \in \{1,\cdots, n\}}^{} {\sigma_{n-1}(E\diagdown \{\tau (x_j)\})} \; x_j\neq x \; \forall x \in E $$ où  $\tau$ est une permutation sur E et que $ \sigma_{n-1} \tau \chi = \sigma  \chi , \; \forall\chi \subseteq E  $. 

Posons  $\mathcal{V}(n) = G$ et $\mathcal{V}'(n) = G'$. Comme $\sigma_{n-1}$  agit sur les parties à n-1 éléments, alors de la représentation vectorielle vue dans le corollaire \ref{c3} on a :

 $$\sigma_{n-1}(\mathcal{V}(n)) \simeq \mathcal{V}'(n)$$ 
 
Or  $$\sigma(\mathcal{V}(n)) \simeq \mathcal{V}'(n)$$ 

Donc $$\sigma_{n-1}(\mathcal{V}(n)) \simeq \sigma(\mathcal{V}(n) $$

Et il existe une permutation $\tau$  sur E telle que :
$$ \sigma(x) = \sigma_{n-1}(\tau(x)) $$

Et comme : $$ \tau(x)= \bigcap_{j \in \{1,\cdots, n\}}^{} {E\diagdown \{\tau (x_j)\}} \; x_j\neq x \;  \forall x \in E $$

Alors : 

$$ \sigma(x) = \sigma_{n-1}(\tau(x))= \bigcap_{j \in \{1,\cdots, n\}}^{} {\sigma_{n-1}(E\diagdown \{\tau (x_j)\})} \; x_j\neq x \;  \forall x \in E $$ 

D'où le résultat.$ \checkmark$

\newpage

 \section{Cas où les graphes  sont multicolores  }
 
Entre autres résultats de la première section, les Théorèmes \ref{t2} et \ref{th mes}  restent  vrais pour les graphes multicolores : Pour leurs démonstrations on se sert - comme  on l'a vu dans la page \pageref{2D1}  - du lemme \ref{l2}.  $\checkmark$

\begin{cor} \label{ms}

 Soit M et N  deux matrices\textbf{ symétriques} de deux transformations linéaires  sur un K-espace vectoriel E de dimension finie $(e_1, e_2, \cdots , e_n)$, $n  \geq 3$.

Soit $M_i$ et $N_i$  les deux matrices obtenues à partir de M et de N en supprimant les $i^{eme}$  lignes et les $i^{eme}$ colonnes.

Si $\forall i \in \{1, \cdots , n\} : M_i= \Sigma_i^{t}N_i \Sigma_i $ où $\Sigma_i$  est la matrice d'une permutation sur  $\{e_1, e_2, \cdots , e_n\} \setminus \{e_i\}$, alors il existe une matrice $\Sigma$ d'une permutation sur  $\{e_1, e_2, \cdots , e_n\} $ telle que $ M=\Sigma^{t} N \Sigma$.
\end{cor}

\textbf{Preuve :} En effet de telles matrices \textbf{symétriques } ne sont que des matrices de graphes multicolores.
$\checkmark$
\newpage

 \section{Généralisation aux multigraphes  multicolores  }
 
 Dans [6], Maurice Pouzet a introduit la notion des multirelations $(R_1 \cdots R_k)$ comme un ensemble de k relations sur une même base E et a généralisé la notion d'isomorphie pour ces multirelations, et par suite les problèmes de reconstruction connus pour les relations se sont généralisés pour les multirelations : En particulier la conjecture d'Ulam pour les multigraphes.
 
 En réponse à la question c [6] de Maurice Pouzet, j'énonce que, entre autres résultats de la première section,  le Théorème \ref{t2}  reste  vrai pour les multigraphes multicolores.
 
 \begin{theo}[Conjecture d'Ulam-Pouzet] [6] 
 
 Soit $(G_1 \cdots G_k)$ et $(G'_1 \cdots G'_k)$ deux multigraphes multicolres sur une même base E de cardinal n au moins égal à 3.
 
 Si $(G_1 \cdots G_k)$ et $(G'_1 \cdots G'_k)$  sont isomorphes sur toute partie de E à n-1 éléments , alors  $(G_1 \cdots G_k)$ et $(G'_1 \cdots G'_k)$ sont isomorphes.
 
 \end{theo}
 
\textbf{ Preuve :}
 
 En effet un multigraphe multicolore $(G_1 \cdots G_k)$  ne sera  qu'un graphe multicolore où chaque arrête est coloré  de k couleurs différentes qu'on considère comme une nouvelle couleur, et le résultat se déduit de la section 2.  $\checkmark$
 
 Exemples de multigraphes multicolores :

\begin{figure}
\centering
\includegraphics[width=100mm,height=100mm,scale=1.5]{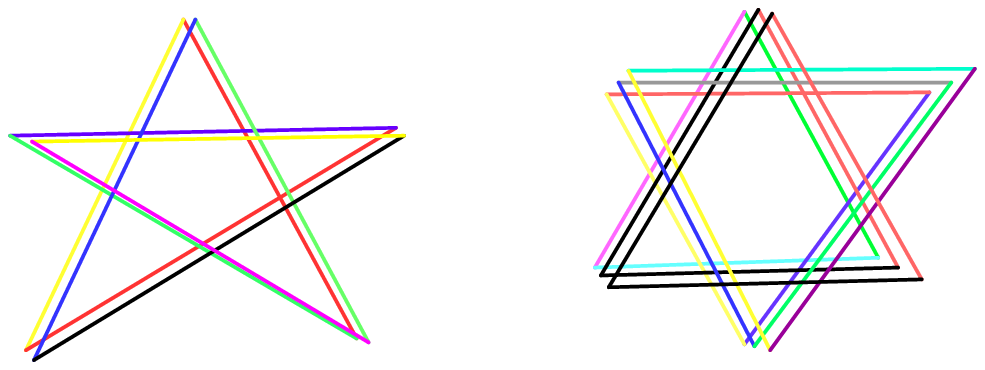}
\caption{Deux multigraphes multicolores}
\label{fig:dm2}
\end{figure}
\newpage

\section{La reconstruction des relations m-aires h-symétriques} 

Dans [6] Maurice Pouzet a définit les relations m-aires héréditairement symétriques  comme suit :

Une relation m-aire sur un ensemble à n éléments ($ n \geq m+1 $)  est héréditairement symétrique (en abrégé h-symétrique) lorsqu'elle prend la même valeur sur deux m-uples définissant la même partie de la base.

Maurice Pouzet a conjecturé que les relations m-aire h-symétriques sur une base à au moins n éléments 
($ n \geq m+1 $) sont Ulam-reconstructibles.

Dans cette section, je donne  une réponse positive à sa conjecture.

\begin{theo}[Conjecture d'Ulam-Pouzet Bis] \label{tup}

Soient G et G' deux relations m-aires h-symétriques sur une même base E de cardinal n au moins égal à $m+1$, $m \geq 2$.

Si G et G' sont (-1)-hypomorphes, alors G et G' sont isomorphes .

\end{theo}

\textbf{Preuve :}

Il suffit de reprendre les mêmes techniques utilisées auparavant  pour les graphes, de prendre les mesures $\mu_i^m$  au lieu de $\mu_i^2$, et de considérer l'égalité :

\begin{eqnarray*}
\bigcup_{i \in \{1,\cdots ,m \} }^{}\sigma(x_i)  =  \bigcap_{k \in \{1,\cdots ,n\} \diagdown \{1,\cdots ,m\}  }^{} \sigma_{n-1}(E\diagdown \{x_k\})
\end{eqnarray*}

Qu'on a utilisé pour $m=2$.

$\checkmark$

\newpage

\section{La limite de validité de l'utilisation des $\infty$ -permutations et l'introduction des dilatations} \label{limite}

Les mesures et les i-permutaions ont permis la reconstruction des multigraphes multicolores \textbf{ finis }à partir de la permutation dont elles sont déduites. Et on peut se demander si l'on peut  généraliser ceci pour le cas des multigraphes multicolores\textbf{ infinis} ? 

La réponse est non : Un contre exemple est donné  par J. Ficher dans [J. F] en construisant deux graphes infinis (-1)-hypomorphes mais non isomorphes.

Voici la raison pour laquelle l'utilisation de mes techniques ne sont pas toujours valables pour les graphes infinis:

D'abord notons $\sigma^*$  une permutation définie sur les parties  propres et maximales de la base  E infinie ( $\sigma^*$ jouera le rôle de $\sigma_{n-1}$ pour le cas fini). $\sigma^*$ est dite une  $\infty$ -permutation.

Notons $\Omega^*_E$ l'ensemble des parties propres et maximales de E .  Et notons  $\Omega_E$ l'ensemble des parties de E abritant au moins une partie propre et maximale de E.
 
Définissons $\Omega_E/ G$  la relation sur l'ensemble $\Omega_E$ des parties de E abritant au moins une partie propre et maximale de E  par :  $ M \equiv N $  mod G   si et seulement si $ G|M  \simeq G|N$. De même on définit $\Omega_E/ G'$ .

Au début, j'ai cherché  de trouver une permutation $\sigma$  dont est déduite la $\infty$ -permutation $\sigma^*$, on peut penser comme j'ai fait dans  la démonstration 2 du Théorème \ref{t1},  de poser :

$$\sigma(x)= E \diagdown \sigma^{*}(E \diagdown \{x\}) \; \forall x \in E $$

Et d'écrire   : $$ \sigma(x_i) =  \bigcap_{x_j \in E \diagdown \{x_i\}}^{} {\sigma^*}(E\diagdown \{x_j\}) \;   $$

Mais de cette égalité,  les graphes infinis trouvés par [J.F] Fisher seront isomorphes, ce qui n'est pas le cas.

\textbf{La raison est simple} : Dans le cas fini, les $\sigma_{n-1}$ font partie des transformations qui envoient toute partie propre et maximale en une autre partie de  sorte  que les mesures reste conservées, mais dans le cas infini, une partie propre et maximale peut être transformée en une partie maximale tout en conservant toutes les mesures (comme dans les contre exemples de [J.F] J.Fisher ).

 \textbf{Physiquement interprété : La force partielle -due à l'explosion-  agissant sur les particules d'une partie propre et maximale a tendance à les répartir dans la base toute entière   tout en conservant toutes les mesures.}

D'où l'idée de considérer  la \textbf{dilatation} $\hat{\sigma^{*}}$   au lieu de  $\sigma^*$ :

La \textbf{dilatation} $\hat{\sigma^{*}}$  de $\sigma^{*}$ définie de  $\Omega^*_E$ sur  $\Omega_E$ par : $\hat{\sigma^{*}}(\chi)= \sigma^{*}(\chi )$ si $\sigma^{*}(\chi) \neq E$  dans  $\Omega_E/ G'$ .  Et $\hat{\sigma^{*}}(\chi)= E$ si $G'|E \simeq G'|\sigma^{*} (\chi)$ dans $\Omega_E/ G'$.

Et l'écriture :  $$ \sigma(x_i) =  \bigcap_{x_j \in E \diagdown \{x_i\}}^{} {\hat{\sigma^*}}(E\diagdown \{x_j\}) \;   $$ ne peut être définie que dans le cas où : $$E \neq \hat{\sigma^{*}}\chi  \; \forall \chi \in \Omega^*_E $$

Donc contrairement au cas fini,  où on a vu dans la preuve du  Théorème \ref{t2} que : $$ \sigma(x_i) =  \bigcap_{j \in \{1,\cdots ,n\} \diagdown \{i\}}^{} {\sigma_{n-1}(E\diagdown \{x_j\})} \;  \forall x_i \in E $$

L'écriture  : $$ \sigma(x_i) =  \bigcap_{x_j \in E \diagdown \{x_i\}}^{} {\hat{\sigma^*}}(E\diagdown \{x_j\}) \;   $$ n'est pas toujours assurée et par la suite on ne peut pas  en général  établir un isomorphisme  entre deux multigraphes multicolores sur une base infinie : En effet c'est l'égalité :

 $$ \sigma(x_i) =  \bigcap_{j \in \{1,\cdots ,n\} \diagdown \{i\}}^{} {\sigma_{n-1}(E\diagdown \{x_j\})} \;  \forall x_i \in E $$  utilisé en particulier dans le lemme \ref{l2} qui a permis d'établir l'isomorphie en écrivant l'égalité : 
 
 \begin{eqnarray*}
 \{\sigma(x_i)\}  \cup \{\sigma(x_j)\} =  \bigcap_{k \in \{1,\cdots ,n\} \diagdown \{i,j\}  }^{} \sigma_{n-1}(E\diagdown \{x_k\})
 \end{eqnarray*}
 
Et en montrant à partir de cette égalité que $\sigma$ préserve la mesure $\mu^{2}_{i}$.

\textbf{Remarque (Surprenante)} \label{remarque} :

Si $\sigma$ ne peut  pas être définie comme dans le cas où on a :$$ \emptyset = E \diagdown \hat{\sigma^*}(E \diagdown \{x\})  $$  pour certaines valeurs x de E.

Dans ce cas : $  G'|E \diagdown \{x\} \simeq G|E \diagdown \{x\}\simeq G'|E$   pour certaines valeurs x de E exactement comme dans  les contre exemples de J. Fisher [J. F] !.

Cette remarque va nous être très utile dans section \ref{Cas des graphes infinis}.

 $\checkmark$
  \newpage

\section{Conjectures   }
La partie précédente et les contres exemples de J. Fisher  permettent de conjecturer :

\begin{cj} \label{cjs1}
Soit G un graghe infini de base E,   si $ \forall x \in E,  G|E \setminus \{x\} \nsim G|E$ alors G est Ulam-reconstructible.

\end{cj}

\begin{cj} \label{cjs2}
Soit G un graghe infini de base E,   si $ \forall (x,y) \in E^2, x \neq y ;  G|E \setminus \{x\} \sim G|E \setminus \{y\}$ alors G est Ulam-reconstructible.

\end{cj}

\begin{cj}   \label{cjs3}
Soit G un graghe infini de base E,   si $ \forall (x,y) \in E^2, x \neq y ;  G|E \setminus \{x\} \nsim G|E \setminus \{y\}$ alors G est Ulam-reconstructible.

\end{cj}

\begin{cj}   \label{cjs4}
Si G est un graghe infini de base E tel que $\forall n \in \mathbb{N}, |E_{n}(G)| < \infty$ où $ E_n(G)= \{ x  \in E \diagup  deg(x) = n \} $, alors G est Ulam-reconstructible.

\end{cj}

\begin{cj}   \label{cjs5}
Les conjectures  \ref{cjs1}, \ref{cjs2} , et   \ref{cjs3}   sont vraies pour les multigraphes multicolores infinis  avec les mêmes hypothèses.

\end{cj}

\begin{cj}   \label{cjs6}
Si G est un graghe infini de base E tel que $\forall x \in E,  \; deg(x) < \infty$ , alors G est Ulam-reconstructible.

\end{cj}

\newpage
\section{Cas des graphes infinis} 

              \label{Cas des graphes infinis}
 
 Le but de cette section est de démontrer les Théorèmes \ref{thcjs1},  \ref{thcjs2} \ref{thcjs3} et \ref{thcjs4} suivants que j'ai conjecturés dans la section  précédente :

Dans la section \ref{limite}, j'ai  introduis et définis les notions de $\infty$-permutation  $\sigma^{*}$ et des dilatations $\hat{\sigma^{*}}$, et surtout j'y ai donné les raisons de leurs introductions. Je les rappelle ici :
  
 Notons $\Omega^*_E$ l'ensemble des parties propres et maximales de E . 
 
 Notons  $\Omega_E$ l'ensemble des parties de E abritant au moins une partie propre et maximale de E.
 
 Définissons $\Omega_E/ G$  la relation sur l'ensemble $\Omega_E$ des parties de E abritant au moins une partie propre et maximale de E  par :  $ M \equiv N $  mod G   si et seulement si $ G|M  \simeq G|N$
  
  De même on définit $\Omega_E/ G'$ .
  
Définition : Une $\infty$-permutation $\sigma^{*}$ sur un ensemble E (fini ou  infini)   est une permutation sur l'ensemble $\Omega^*_E$ .

 Pour mieux comprendre, si $\sigma$ est un isomorphisme  entre deux graphes G et G' (-1)-hypomorphes sur une base E, alors on peut poser :
 
 $$ \sigma(x)= E \diagdown \sigma^{*}(E \diagdown \{x\}) \; \forall x \in E   $$
 
 Si $G|E \diagdown \{x\} \simeq G|E$ , alors  $G'|\sigma(E \diagdown \{x \}) \simeq G|E \diagdown \{x\} \simeq G|E \simeq G'|\sigma(E)$ 
 
 Donc :  $$G'|\sigma(E \diagdown \{x \}) \simeq G'|E$$
 
 Or $\sigma(E \diagdown \{x \}) = \sigma^{*}(E \diagdown \{x\})$, donc : $G'|\sigma^{*}(E \diagdown \{x\})  \simeq G'|E$.

 Ceci permet de définir la \textbf{dilatation} $\hat{\sigma^{*}}$  de $\sigma^{*}$ définie de  $\Omega^*_E$ sur  $\Omega_E$ par : $\hat{\sigma^{*}}(\chi)= \sigma^{*}(\chi )$ si $\sigma^{*}(\chi) \neq E$  dans  $\Omega_E/ G'$ .  Et $\hat{\sigma^{*}}(\chi)= E$ si $G'|E \simeq G'|\sigma^{*} (\chi)$ dans $\Omega_E/ G'$.

On a donc :

$$\sigma \chi \simeq \hat{\sigma^{*}}\chi  \; \forall \chi \in \Omega^*_E $$

D'où la question : Si $\hat{\sigma^{*}}$ existe et préserve les mesures $\mu^2_i$, peut-on trouver une permutation $\sigma$ sur E préservant les mesures  $\mu^2_i$ et telle que :
 
 $$\sigma \chi \simeq \hat{\sigma^{*}}\chi  \; \forall \chi \in \Omega^*_E $$

On verra que la réponse est positive dans le cas où : $\hat{\sigma^{*}} = \sigma^{*}$.

\textbf{Remarques : }
  
1- Si E est fini alors $\hat{\sigma^{*}} = \sigma^{*}$.

2- L'idée de l'introduction de la \textbf{dilatation} $\hat{\sigma^{*}}$   est née du fait que dans les contre-exemples de J. Fisher [J F] on a :  $ G|E \diagdown \{x\} \simeq  G'|E \diagdown \{x\} \simeq G'|E $, soit $\hat{\sigma^{*}}(E \diagdown \{x\}) = E $.
 
 $\hat{\sigma^{*}}$ existe et est bien définie si  $\sigma^{*}$  si existe (C'est le cas  si G et G'  sont (-1)-hypomorphes).

 Si $$ \emptyset \neq  E \diagdown \hat{\sigma^{*}}(E \diagdown \{x\}) \; \forall x \in E $$  Alors : 
 
 $\sigma(x)= E \diagdown \hat{\sigma^{*}}(E \diagdown \{x\}) $  sera bien définie $ \forall x \in E $ , et les mêmes techniques  utilisées dans le cas fini montrent que G et G' sont isomorphes.

 Par contre si :$$ \emptyset =  E \diagdown \hat{\sigma^{*}}(E \diagdown \{x\}) \;  $$  pour certaines valeurs x de E, alors comme on l'a vu précédemment on ne peut pas  toujours établir un isomorphisme entre G et G' mais :
 
 $$ \exists x \in E \; \diagup  G|E \diagdown \{x\} \simeq  G'|E $$ 
 
 Car : $$ \emptyset =  E \diagdown \hat{\sigma^{*}}(E \diagdown \{x\}) \; \text{pour certains x de E}  $$ 
 
 Or  : $$ G|E \diagdown \{x\} \simeq  G'|E \diagdown \{x\} $$
 
 Donc : $$ \exists x \in E \; \diagup G'|E \diagdown \{x\} \simeq G'|E $$
 
 De même :  $$ \exists y \in E \; \diagup G|E \diagdown \{y\} \simeq G|E $$

 De ceci on déduit :

 \begin{theo} \label{thcjs1}
 Soit G un graghe infini de base E,   si $ \forall x \in E,  G|E \setminus \{x\} \nsim G|E$ alors G est Ulam-reconstructible.

 \end{theo}

 \begin{theo} \label{thcjs2}
 Soit G un graghe infini de base E,   si $ \forall (x,y) \in E^2, x \neq y ;  G|E \setminus \{x\} \sim G|E \setminus \{y\}$ alors G est Ulam-reconstructible.
 
  \end{theo}
\begin{theo}   \label{thcjs3}
Soit G un graghe infini de base E,   si $ \forall (x,y) \in E^2, x \neq y ;  G|E \setminus \{x\} \nsim G|E \setminus \{y\}$ alors G est Ulam-reconstructible.

\end{theo}  

\textbf{Preuve : }

Si G' est (-1)-hypomorphe mais non isomorphe à G, alors :  $$\exists x \in E  \;\diagup \; \emptyset =  E \diagdown \hat{\sigma^{*}}(E \diagdown \{x\}) \;  $$ 

Soit : $$ G|E \setminus \{x\} \sim G'|E \setminus \{x\} \sim G|E$$

Et comme :  $ \forall (x,y) \in E^2, x \neq y ;  G|E \setminus \{x\} \nsim G|E \setminus \{y\}$, alors :

 $ \forall (x,y) \in E^2, x \neq y ;  G'|E \setminus \{x\} \nsim G'|E \setminus \{y\}$
 
 On en déduit que :  $$\forall y \in E \diagdown \{x\}, \; \emptyset \neq  E \diagdown \hat{\sigma^{*}}(E \diagdown \{y\}) \;  $$  
 
Que : $$ \sigma^{*}(E \diagdown \{x\})= E \diagdown \{x\}$$
 
Et que $ \sigma(z)=E  \diagdown \sigma^{*}(E \diagdown \{z\})  $ est bien définie  pour tout élément z de E. Et par suite,  les mêmes techniques  utilisées dans le cas fini montrent que G et G' sont isomorphes.
 
  \begin{theo}   \label{thcjs4}
  Soit G un graghe infini de base E tel que $\forall n \in \mathbb{N}, |E_{n}(G)| < \infty$ où $ E_n(G)= \{ x  \in E \diagup  deg(x) = n \} $, alors G est Ulam-reconstructible.

  \end{theo}
  
\textbf{  Preuve : }

Se déduit du Théorème \ref{thcjs1}

  $\checkmark$
  \newpage
  \section{Les graphes infinis et localement finis }\label{partie 7}
  
  Le but de cette section est de démontrer le Théorème \ref{thcjs6} suivant.

  \begin{theo}   \label{thcjs6}
  Si G est un graghe infini de base E tel que $\forall x \in E, \;  deg(x) < \infty$ , alors G est Ulam-reconstructible.

  \end{theo}
  
  \textbf{Preuve :}
  
  On a vu  dans le Théorème \ref{thcjs1} que si G' est (-1)-hypomorphe mais non isomorphe à G alors :
  
  $$ \exists x \in E \; \diagup G'|E \diagdown \{x\} \simeq G'|E $$
  Posons $E=F_0$.
  
  Posons $x=x_1$ et $E \diagdown \{x_1\}=F_1$.
  
  On a donc : $$ G'|F_0 \simeq G'|F_1 $$
  
  Or en changeant la base E de G par $F_1$, comme G' et G sont  (-1)-hypomorphes sur $F_1$ mais non isomorphes, alors :
  $$ \exists x_2 \in F_1 \; \diagup G'|F_1 \diagdown \{x_2\} \simeq G'|F_1$$
  
  En continuant ce procédé; il existe donc une suite décroissante  d'ensembles $F_i$ tels que : 
   $$  G'|F_i \simeq G'|F_{i+1} , F_i \supsetneq F_{i+1},  |F_i \diagdown F_{i+1}|=1 $$

  Or  par (-1)-hypomorphie et par construction  il existe une suite  décroissante  d'ensembles $E_i$ tels que : 
     $$ G|E_{i+1} \simeq G|E_i \simeq G'|F_i , E_i \supsetneq E_{i+1} , |E_i  \diagdown E_{i+1}|=1 , \forall i \in \mathbb{N^*}  \text{avec} \; G|E_i \nsim G|E $$

On a donc :
 
 $ G|E_{i+1} \simeq G|E_i \simeq G'|F_i , E_i \supsetneq E_{i+1} , |E_i \diagdown E_{i+1}|=1 , 
 \forall i \in \mathbb{N^*} ,  \exists x_i \in E_i \; \diagup \; G|E_i \diagdown \{x_i\} \nsim G|E $

%$E \diagdown E_i$ est une suite d'ensembles croissante  pour l'inclusion car  $E_i$  est décroissante. De plus $E \diagdown E_i$  est majorée par E.
%
%Donc  $E \diagdown E_i$ a une limite L. L ne peut pas être égal à E car sinon  $E_i$ aura $\emptyset$  pour limite et $\text{card}E_i$  sera fini pour i assez grand, ce qui est impossible  car : $$ G|E_{i+1} \simeq G|E_i \simeq G'|F_0 , \forall i \in \mathbb{N^*}  $$

Dans l'ensemble des parties $\chi \subset  E$ telles que $ G|\chi \simeq G|E_1 $, les $E_i$  ont une limite $E_m $. (En fait, on peut voir facilement que le dit ensemble coïncide  avec son adhérence pour l'inclusion).

Par minimalité, $  G|E_m \diagdown \{x\} \simeq G|E, \, \forall x \in E_m $.

\textbf{$1^{er}$ Cas} : Si tout point de E est au moins de degré un  relativement à G.

De l'isomorphisme $  G|E_m \diagdown \{x\} \simeq G|E, \, \forall x \in E_m $,  on déduit qu' il existe un point x de E de degré nul ce qui est absurde.

\textbf{$2^{eme}$ Cas} :  Si il existe au moins un point de E de degré nul  relativement à G.

- Si E a un nombre infini de points de degré nul relativement à G alors $ E_m $ a un nombre infini de points de degré nul relativement à G, en effet : comme  $  G|E_m \diagdown \{x\} \simeq G|E, \; \forall x \in E_m $, alors $  G|E_m \diagdown \{x\}$  a un nombre infini de points $x_i$ de degré nul et pour un nombre infini de ces $x_i$  on a : $G(x,x_i)=0$, car sinon le point x sera de\textbf{ degré infini}, ce qui est absurde, donc on a un nombre infini de points de degré nul simultanément pour  $  G|E_m$ et pour  $  G|E$. 

Comme $  G|E_m \diagdown \{x\} \simeq G|E, \, \forall x \in E_m $, alors pour x de degré nul on aura : $  G|E_m  \simeq G|E $, et par suite :  $  G|E \diagdown \{x\} \simeq G|E, \, \forall x \in E $, soit $ \forall (x,y) \in E^2, x \neq y ;  G|E \setminus \{x\} \sim G|E \setminus \{y\}$ .

Or dans ce cas G est Ulam-reconstructible d'après le Théorème \ref{thcjs2}. Soit donc G et G' sont isomorphes, ce qui est absurde.

- Si  E a un nombre fini de points de degré nul relativement à G :

* Ou bien  tout point de E est au moins de degré un  relativement à G' on aboutit à une absurdité comme dans le $1^{er}$ Cas.

*  Ou bien E a un nombre infini de points de degré nul relativement à G' on aboutit à une absurdité comme ci-dessus.

*  Ou bien E a un nombre fini de points de degré nul relativement à G', par construction comme $  G'|F_i \simeq G'|F_{i+1} \simeq G'|E , \forall i  \in \mathbb{N}^{}$ , alors chaque $F_i$  a  des points de degré nul relativement à G', et par isomorphie : $G|E_{i+1} \simeq G|E_i \simeq G'|F_i $, il en sera de même de chaque  $E_i$ relativement à G, et de $E_m$ relativement à G. Or comme $ \forall (x,y) \in E_m^2,  x \neq y$ on a  $G|E_m \setminus \{x\} \sim G|E_m \setminus \{y\}$, alors  $ E_m $  a un nombre infini de points de degré nul relativement à G, et on aboutit à une absurdité comme ci-dessus.

$\checkmark $

\newpage

\section{Vision dynamique de la reconstruction des graphes et conclusion} 

\textbf{Commentaire sur le lemme } \ref{l2}: \label{VD}

Si G est un graphe sur une base E de cardinal $n \geq 3$. J'ai considéré la restriction de G à toute partie F de E comme une mesure sur F.

Si G' est un autre graphe sur E. Et si $\sigma $ est un isomorphisme entre G et G' sur E, alors il existe une permutation $\sigma^*$ (ou $\sigma_{n-1}$  si E est de cardinal n)  sur les parties propres et maximales et préservant les mesures.

On a définit la \textbf{dilatation} $\hat{\sigma^{*}}$  de $\sigma^{*}$ définie de  $\Omega^*_E$ sur  $\Omega_E$ par : $\hat{\sigma^{*}}(\chi)= \sigma^{*}(\chi )$ si $\sigma^{*}(\chi) \neq E$  dans  $\Omega_E/ G'$ .  Et $\hat{\sigma^{*}}(\chi)= E$ si $G'|E \simeq G'|\sigma^{*} (\chi)$ dans $\Omega_E/ G'$.  où on a noté : 

  $\Omega^*_E$ l'ensemble des parties propres et maximales de E . 
 
  $\Omega_E$ l'ensemble des parties de E abritant au moins une partie propre et maximale de E.

On a vu que si E est de cardinal fini, alors   $\sigma^*=\hat{\sigma^{*}}$

\textbf{Le problème inverse }: Si l'action $\hat{\sigma^{*}}$  transforme une  partie propre et maximale en une autre partie contenant une partie propre et maximale et en préservant  les mesures $\mu^2_i$. Peut-on trouver une permutation $\sigma $  sur E qui conserve les mesures $\mu^2_i$ entre G et G'?

\underline{Dans le cas où E est fini :}

Le lemme \ref{l2} affirme l'existence de $\sigma\tau$  qui conserve les mesures $\mu^2_i$ entre G et G' avec $\sigma^*$  déduite de $\sigma$ .

On peut interpréter cela comme ceci : Si les points $x_i$ sont des particules, sous l'action des $\sigma^*(E\diagdown \{x_j\})$ , les points $x_i$ vont se positionner aux points $\sigma(x_i)$ de sorte que les mesures $\mu^2_i$ entre G et G' restent conservées.

\underline{Dans le cas où E est infini :}

On n'a pu affirmer l'existence d'un isomorphisme  $\sigma$ entre G et G' que dans le cas où aucune partie propre et maximale $E\diagdown \{x_j\}$ n'est envoyée sur E  sous l'action de $\hat{\sigma^{*}}$ - En particulier  si la restriction de G' à E n'est pas  isomorphe à la restriction de G' sur une partie propre et maximale, ce qui n'est  pas le cas dans les contre-exemples de J. Fisher - .

Je rappelle que  si la restriction de G' à E n'est pas  isomorphe à la restriction de G' sur une partie propre et maximale, alors $\sigma^*=\hat{\sigma^{*}}$ comme dans le cas fini.

\textbf{Vision dynamique de le reconstruction des graphes :}

Dans les cas où on a  pu démontrer la conjecture d'Ulam , on a : $\sigma^*=\hat{\sigma^{*}}$.

Si les points $x_i$ sont des particules, sous l'action des $\sigma^*(E\diagdown \{x_j\})$ , les points $x_i$ vont se positionner aux points $\sigma(x_i)$ de sorte que les mesures $\mu^2_i$ entre G et G' restent conservées.

On peut donc voir $\sigma$  comme l'action  résultante  des actions $\sigma^*(E\diagdown \{x_j\})$  où $x_j \in E$.

 \newpage

\section*{ Référence bibliographique }
 \addcontentsline{toc}{section}{ Référence bibliographique }

[1] Bondy,J.A., and R.L. Hemminger

Graph reconstruction, a survey. J. Graph Theory1 (1977), 227-268

[2] P. J. Cameron

Stories from the age of reconstruction, Congr. Numer. 113 (1996) 31-41. Festschrift for C. St. J.A. Nash-Williams

[J. Ficher]

A counterexample to the countable version of a conjecture of Ulam. Journal of combinatorial theory 7, 364-365 (1969).

[3] F. Haray, E. Palmer

On the problem of the reconstruction of a tournament from subtournaments. Mh. Math. 71 (1967), 14-23.

[4] Kelly P. J.

A congruence theorem for trees, Pacific J. Math. 7 (1957) 961-968

[5] V. B. Mnukhin.

The k-orbit reconstruction and the orbit algebra, Acta Appl. Math. 29(1-2).(1992) 83-117.Interactions between algebra and combinatorics.

[6] M.Pouzet

Application d’une propriété combinatoire des parties d’un ensemble aux groupes et aux relations.Math. Zeitschrirt. 150 (1976). P. 117-134

[7] M. Pouzet

Relations non reconstructible par leurs restrictions, Journal of combinatorial Theory, Series B 26, 22-34(1979)

[8]M. Pouzet, N. M. Thiéry

Invariants algébriques de graphes et reconstruction, C. R. Acad. Sci. Paris, t. 333, Série I, p. 821-826, 2001.

[9]P. K. Stockmeyer 

The falsity of the reconstruction conjecture for tournaments. J. Graph Theory.1. 1977. p.19-25.

[10]S. M. Ulam  

“ A Collection of Mathematical Problems,”Interscience, New York, 1960

\begin{flushright}
Mohamed sghiar

msghiar21@gmail.com

Tel : 0033(0)953163155\& 0033(0)669753590.

\end{flushright}

\newpage

\end{document}